\documentclass[12pt]{amsart}
\usepackage{amsmath,amssymb,amsthm}

\newtheorem{theorem}{Theorem}[section]
\newtheorem{proposition}[theorem]{Proposition}
\newtheorem{lemma}[theorem]{Lemma}
\newtheorem{corollary}[theorem]{Corollary}
\newtheorem{remark}[theorem]{Remark}
\newtheorem{definition}[theorem]{Definition}
\newtheorem{example}[theorem]{Example}

\newcommand{\C}{{\mathbb C}}

\newcommand{\Q}{{\mathbb Q}}
\newcommand{\R}{{\mathbb R}}
\newcommand{\Z}{{\mathbb Z}}

\newcommand{\matO}{{\mathcal O}}

\newcommand{\too}{\longrightarrow}

\renewcommand{\L}{{\mathcal L}}

\DeclareMathOperator{\length}{length}

\newcommand{\Pl}{{\mathbb P}}
\newcommand{\Qbar}{\bar{\Q}}

\DeclareMathOperator{\Spec}{Spec}
\DeclareMathOperator{\rad}{rad}\DeclareMathOperator{\Tor}{Tor}

 \DeclareMathOperator{\Pic}{Pic}

\DeclareMathOperator{\proj}{proj} 
 
\DeclareMathOperator{\Div}{div} \DeclareMathOperator{\Ar}{Ar}
\DeclareMathOperator{\finite}{finite} \DeclareMathOperator{\N}{N}

\newcommand{\vfi}{{\varphi}}

\newcommand{\bC}{{\mathbb C}}
\newcommand{\bR}{{\mathbb R}}
\newcommand{\bQ}{{\mathbb Q}}

\newcommand{\cO}{\mathcal{O}}

\newfont{\cyrr}{wncyr10}
\title[n-dimensional Mahler formula]{Mahler formula for self maps on the
n-dimensional projective space}
\author{J.~Pineiro}
\thanks{I'm very grateful to Professor Lucien Szpiro for making me aware,
between many other ideas, of the geometric aspects of the Mahler
formula.  I would like also to thank Michael Tepper for pointing out
some typos in the original version of this manuscript}


\address{
Department of Mathematics and Computer Science\\
Bronx Community College\\
University Ave. and West 181 Street\\
Bronx, NY 10453}
\email{jorge.pineiro@bcc.cuny.edu}

\begin{document}

\begin{abstract}
The Mahler formula gives an expression for the height of an
algebraic number, as the integral of the log of the minimal equation
with respect to the Haar measure on the circle.  In the present work
we prove that a similar result holds for nice self maps on the
n-dimensional projective space.  The height of the number is
replaced by the canonical height of a hypersurface, and the Haar
measure is replaced by the canonical invariant measure.

\end{abstract}

\maketitle
\section{Introduction}
A classical formula of Mahler \cite{mahler} states that, if $P=(\lambda:1)
\neq \infty$ is a point
in $\Pl^1_{\bar{\Q}}$, the naive height $h_{nv}(P)$
can be related to the integral of the log of the minimal equation
$F$ of $\lambda$, with respect to the Haar measure $d\theta$ on
$S^1$.  The identity we find in this case is:
$$h_{nv}(P)=\frac{1}{\deg(F)} \int_{S^1} \log|F(z)|d\theta.$$ Now, the
naive height is the canonical height (c.f. later) associated to the
morphisms $\phi_n : t \rightarrow t^n$ on $\Pl^1$. This means that
$h_{nv}(\phi_n(t))=nh_{nv}(t)$ and $h_{nv}(t) \geq 0$. On the other
hand the Haar measure $d\theta$ on $S^1$ is invariant under the
action of this endomorphisms, in the sense that $\phi_n^* d\theta =
n d\theta$ and $\phi_{n*} d\theta = d\theta$. It was already
established in \cite{sz-t-p} that a more general equation is true
for arbitrary dynamical systems on $\Pl^1$.  The terms ``bad
reduction'' and ``integral at a finite place'' appear for the first
time while dealing with dynamical systems in dimension one.  We
prove in the following work that the situation is similar for
dynamical systems on $\Pl^n$, provided that the morphism admits a
``good model''.  Suppose that we are working with a number field
$K$. The general Mahler formula we prove states that, if $\vfi :
\Pl_K^n \rightarrow \Pl_K^n$ admits a model $\vfi=(p_0:...:p_n) :
\Pl_K^n \rightarrow \Pl_K^n$, such that $(p_0,...,p_n)$ represents a
regular sequence inside $\cO_K[T_0,...,T_n]$ then:
$$h_{\vfi}(D^+)-h_{\vfi}(D^-)=\sum_{v/\infty} \int_{\Pl^n_{\C}} \log
|F|_v d\mu_{\varphi,v} + E(F,v \finite) $$ where we have the
following:
\begin{enumerate}
\item $F=F^+/F^-$ is a rational function on $\Pl^n_K$.
\item $h_{\vfi}(D^+)$ and $h_{\vfi}(D^-)$ represent the canonical heights of
the
cycles $D^+$ and $D^-$ respectively, where $\Div(F)=D^+-D^-$.
\item For every place $v$ of $K$ at infinity, $d\mu_{\vfi,v}$ represents an
invariant measure relative to $\vfi$ on $\Pl^n_{\C}$.
\item Even though the sequence $(p_0,...,p_n)$ is regular in
$\cO_K[T_0,...,T_n]$, the map $\vfi=(p_0:...:p_n) : \Pl^n_K
\rightarrow \Pl^n_K$ may not be a well defined map on
$\Pl^n_{\cO_K}$.  For example $\vfi : \Pl_{\Q}^2 \rightarrow
\Pl^2_{\Q}$, defined over $\Q$ as
$\vfi(x:y:z)=(y^2-3z^2:x^2-3y^2:zy)$ does not extend to a map on
$\Pl^2_{\bar{\Q}_3}$.  The term $E(F,v \finite)$ in the formula, is
arising from the blow-up we may need to do in order to extend the
map $\vfi$ to an integral model. It depends in fact only on a finite
number of finite places, which we call places of bad reduction. When
the map $\vfi : \Pl^n_K \rightarrow \Pl^n_K$ defines a map (which we
are calling here with the same name) $\vfi : \Pl^n_{\cO_K}
\rightarrow \Pl^n_{\cO_K}$, the term $E(F,v \finite)=0$.
\end{enumerate}
We will explain several particular cases and consequences of our formula;
the classical
formula of Mahler will be between them.  Let us assume for the moment that
our map admits a model such
that $h_{\vfi}(D^-)=0$.  This condition will prove to be natural for
polynomial functions in dimension one.
If we pick the equation $F$ such that $v(F)=0$ for every finite
place $v$ (the valuation $v$ naturally extends to rational functions
on $\Pl^n_K$), and $E(F,v \finite) < 0$, we get the inequality:
$$h_{\vfi}(D^+) \leq \sum_{v/\infty} \int_{\Pl^n_{\C}} \log
|F|_v d\mu_{\varphi,v}.$$ When the map $\vfi$ has good reduction
everywhere ($E(F,v \finite) = 0$) and $v(F)=0$ for every finite $v$,
the above inequality becomes an equality.  Particular cases of this
formula can be found in \cite{maillot}.\\
Let's see now the dimension one case, which was treated in
\cite{sz-t-p}.  Suppose that $\vfi=(p_0:p_1): \Pl^1 \rightarrow
\Pl^1$ is a map on the Riemann sphere and $F$ is a polynomial
equation.  We can always change coordinates to get $T_1/p_1$, which
will make $h_{\vfi}(D^-)=h_{\vfi}(\infty)=0$.  Also by base change
we can assume that $(p_0,p_1)$ is a regular sequence.  As a
consequence of Proposition \ref{convergence of finite place int} and
definition \ref{definition of int at finite place} we will be able
to consider the term $E(F,v \finite)$ as sum of ``integrals'' over
the finite places of $K$, i.e. $E(F,v \finite)=\sum_{v
\finite}\int_{\Pl^n_{\C_v}} \log |F| d\mu_{\vfi,v}$.  The measure
$d\mu_{\vfi}$ at each place over infinity is nothing but the Brolin
measure \cite{Brolin}, further studied by other authors such as,
Lyubich \cite{Lyubich} and Freire, Lopez, and Ma\~{n}e \cite{Mane}.
Taking a point $P=(\lambda:1) \neq \infty$ in $\Pl_{\bar{K}}^1$ and
$F$ the minimal equation of $\lambda$ over $\cO_K$ (which we assume
with no common factors, i.e. $v(F)=0$ at every finite valuation),
the formula we found takes a more symmetric shape,
$$h_{\vfi}(D^+) = \sum_{v } \int_{\Pl^1_{\C}} \log
|F|_v d\mu_{\varphi,v}.$$
When the map
$\vfi=(p_0:p_1)$ is defined over $\Q$ and has good reduction
everywhere (which is the case of the morphisms $\phi_n : t
\rightarrow t^n$), we can take $F$ to be the minimal equation of
$\lambda \in \bar{\Q}$ over $\Z$ and get:
$$h_{\vfi}(P)=\frac{1}{\deg(F)}  \int_{\Pl^1_\C} \log |F|
d\mu_{\varphi}.$$ This last expression is the Mahler formula for
rational morphisms on $\Pl^1$ in absence of bad reduction.\\
Going back to the general case, in this paper we work with a map
$\vfi=(p_0:...:p_n) : \Pl^n_K \rightarrow \Pl^n_K$. It may not be
possible to extend this map to a well defined map on
$\Pl^n_{\cO_K}$, exactly because the $p_i$ may have common zeroes
along a subscheme $Y_1$ of $\Pl^n_K$. The technique we will use is
arithmetic intersection theory. We will work with arithmetic
varieties $X_k$ and rational maps $\sigma_k : X_k \rightarrow
\Pl^n_K$, determined by blowing-up subschemes $Y_k$ in $\Pl^n_K$. We
will establish the equality of cycles:
\begin{equation*} \label{inteF_k}
\begin{split}
\Div(F_k) =  D_k  - \deg(F) \infty_k & + \sum_{v,i} x_{v,i,k}
C_{v,i,k}  \\ & -\deg(F) \sum_{v,i} y_{v,i,k} C_{v,i,k} +
\sum_{\finite \,v} v(F) X_{v,k}.
\end{split}
\end{equation*}
where
\begin{enumerate}
\item $F$ is a polynomial in the n variables variables
$T_0/T_n,...,T_{n-1}/T_n$ and $F_k=\sigma^*_k F$;
\item $\Div(F)=D-\deg(F)\infty + \sum_{\finite \,v} v(F) X_v$ and $D_k$ is the proper transform of
$D$ by $\sigma_k$;
\item the divisor $\infty_k$ is determined by the equation
$\sigma_k^*T_n=0$ in $X_k$ and;
\item the
$C_{v,i,k}$ are the different components of the exceptional divisor
of the blow up.
\end{enumerate}
Then we will intersect both sides of the above expression with a
particular class of curves on $X_k$ and finally we will get the
Mahler formula from a limit argument on $k$.
\subsection{Notation and conventions} Unless otherwise stated $K$ will
denote a number field with ring of integers $\cO_K$.  For a place
$v$ of $K$, $\C_v=\bar{K}_v$ will denote the completion of the
localization ${K_v}$ at $v$.  $\Pl_K^n$ will denote the
n-dimensional projective space over $K$ and similar for
$\Pl_{\cO_K}^n$.  The symbol $\L$ will be used to denote line
bundles on different kind of varieties. $V$ will denote an algebraic
variety (must of the time projective) of dimension n. $M$ will
denote a complex projective variety of dimension n and, if $\L$ is a
line bundle on $M$, the term $c_1(\L,\|.\|)=c_1(\bar{\L})$ will
denote a $(1,1)$ current, similar to the first Chern form of $\L$.
In the presentation of the arithmetic as well as the geometric
interception theory, $X$ will denote a Macaulay arithmetic variety
of absolute dimension $n+1$ over $\Spec(\cO_K)$.  This means that
all local rings are Macaulay and there exist a flat, proper and
finite type map $f: X \rightarrow \Spec(\cO_K)$ whose fibres $X_{v}$
over the places $v$ of $K$ are projective varieties of dimension
$n$. For an arithmetic variety $X$ and  a line bundle $\L$ on $X$,
we denote by $\L_v=\L \otimes_{\cO(X)} \Spec(K_v)$ the restriction
of $\L$ to the fibre $X_v$. The line bundle $\L$ will comes
sometimes equipped with hermitian metrics $\|.\|_{P,v}$ on the
fibres $\L_{P,v}$ over each point $P \in X_v$. For a set
$\L_1,...,\L_i$ of metrized line bundles on an arithmetic variety
$X$, the expression
$\hat{\deg}(\hat{c}_1(\L_1)...\hat{c}_1(\L_i)|Z)$ will represent the
arithmetic intersection degree of the line bundles $\L_1,...,\L_i$
over a cycle $Z \subset X$ of dimension $i$. We will pay special
attention, in sections 4, 5 and 6 of this paper, to the arithmetic
variety $X=\Pl^n_{\cO_K}$ as well as arithmetic varieties $X_k$ that
arise from blowing-up subschemes $Y_k$ of $X$. By a model for a map
$\vfi : \Pl_K^n \rightarrow \Pl^n_K$ we mean a system
$(p_0:...:p_n)$ of polynomials in $K[T_0,...,T_n]$ representing the
map in the coordinates $(T_0:...:T_n)$.  The ideal generated by a
system of polynomials $p_0,...,p_m$ will be denoted by $\langle
p_0,...,p_m \rangle$ and the symbol $\rad(I)$ will be used to denote
the radical of any ideal $I$.

\section{Self maps on algebraic varieties}
In this section we will give some examples of self maps on algebraic
varieties.  Let $\varphi: V \rightarrow V$ be a map of the algebraic
variety $V$ to itself. Under certain conditions on the variety $V$
(existence of a line bundle $\L$ with good properties), we will
associate to $V$ and $\L$ a canonical height function and a
canonical measure. The canonical height will be a generalization of
both, the Neron-Tate on abelian varieties and the naive height on
$\Pl^1$.  The canonical measure will be a generalization of Brolin's
measure for maps on $\Pl^1$.  Let's start by giving some examples of
self maps on algebraic varieties.
\begin{example} \label{n-projspace}
Suppose that $K$ is a field. A map $\varphi : \Pl_K^n \rightarrow
\Pl_K^n$ of algebraic degree $d$ is given by a set of degree $d$
homogenous polynomials $p_0(T_0,...,T_n),...,p_n(T_0,...,T_n) \in
K[T_0,...,T_n]$, such that $\rad (\langle p_0,...,p_n
\rangle)=\langle T_0,..., T_n \rangle$.
\end{example}

\begin{example}
Let $A$ be an abelian variety. The multiplication by n morphism $[n]
: A \rightarrow A$ represents a self map on $A$.
\end{example}

\begin{example}
Consider the smooth toric variety $\Pl(\Delta)$ defined over
$\Qbar$.  For each $p \geq 2$ there is a morphism $[p]: \Pl(\Delta)
\rightarrow \Pl(\Delta)$ associated with the multiplication by $p$
in $\Delta$.
\end{example}

\subsection{Canonical height}
Suppose that $V$ is a projective variety defined over a number field
$K$, and $\varphi : V \too V$ a morphism, with the property that
there exist a line bundle $\L$ on $V$ and a real number $\alpha>1$
such that $\varphi^*\L \xrightarrow{\sim} \L^{\alpha}$.  Suppose
also that $h_{\L}$ represents a height function associated to $\L$
(for a detailed discussion see B-3 in \cite{intro-diophantine}).
Then we can find (see for example B-4 in \cite{intro-diophantine}) a
positive height function $h_{\varphi}$ on $V(\bar{K})$, defined as
the limit $h_\varphi(P) = \lim_{k\rightarrow\infty}\frac{
h_{\L}(\varphi^k(P))}{\alpha^k}$ with the properties:
\begin{enumerate}
  \item $h_\varphi$ satisfies Northcott's theorem: points with coordinates
in
  $\bar{K}$ with bounded degree and bounded height are finite in number.

   \item $h_\varphi(\varphi(P)) = \alpha h_\varphi(P).$

   \item $h_\varphi$ is a non-negative function.

   \item $h_\varphi(P) = 0$ if and only if $P$ has a finite forward
     orbit under iteration of the map $\varphi$.

   \item $|h_\varphi(P) - h_{\L}(P)|$ is bounded on
   $V(\bar{\Q})$.
\end{enumerate}

Condition (iv) above expresses the general ideal that canonical
heights should reflect the complexity of algebraic cycles under
iteration of maps. With the use of the arithmetic intersection
theory \cite{zhangadelic} we can extend the positive function
$h_{\vfi}$ to algebraic cycles in $V$.  The value $h_{\vfi}(Y)$ will
denote the canonical height of a $p-$cycle $Y$ inside $V$ and we
will have some similar properties (Check \cite{zhangadelic} theorem
2.4):
\begin{enumerate}
\item $h_{\vfi}(Y) \geq 0$.
\item $h_{\vfi}(Y)$ satisfy the functional equation
$h_{\vfi}(\vfi_*Y)= \alpha h_{\vfi}(Y)$.
\item If the orbit $\{Y,f(Y),...\}$ is finite, then
$h_{\vfi}(Y)=0$.
\end{enumerate}

In section 3 of this paper, we will develop the necessary arithmetic
intersection theory to define the height $h_{\L_1,...,\L_n}$
attached to a collection of metrized line bundles $(\L_i,\|.\|_i)$.
In \cite{zhangadelic} this results are generalized to define heights
associated to limits of metrics, the so-called adelic metrized line
bundles.  The canonical height will naturally arise associated to a
special kind of adelic metric on our line bundle $\L$, which will be
called the canonical metric.

\subsection{Canonical metrics}Let $K$ be an number field and $v$ a place of
$K$.
  Consider the projective variety $V$ defined over $\bar{K}_v$ and $\L$ a
line bundle on $V$ such that $\phi : \L^{\alpha}
\xrightarrow{\sim} \varphi^*\L$ for some $\alpha
>1$. Assume that we have chosen a continuous and bounded metric $\|.\|_v$ on
each fibre of $\L_v$.  The following theorem is
due to Shouwu Zhang \cite{zhangadelic}:
\begin{theorem} \label{canonical metric}
The sequence defined recurrently by $\|.\|_{v,1}=\|.\|_v$ and
$\|.\|_{v,n}=(\phi^* \varphi^* \|.\|_{v,n-1})^{1/\alpha}$ for $n
>1$, converge uniformly on $V(\bar{K}_v)$ to a metric
$\|.\|_{\varphi,v}$ (independent of the choice of $\|.\|_{v,1}$) on
$\L_v$ which satisfies the equation
$\|.\|_{\varphi,v}=(\phi^*\varphi^*\|.\|_{\varphi,v})^{1/\alpha}$.
\end{theorem}
\begin{proof}
See theorem (2.2) in S. Zhang \cite{zhangadelic}. Denote by $h$ the
continuous function $ \log \frac{\|.\|_2}{\|.\|_1}$ on
$V(\bar{K}_v)$. Then
$$\log\|.\|_n=\log\|.\|_1 + \sum_{k=0}^{n-2}(\frac{1}{\alpha}
\phi^*\varphi^*)^k h.$$ Since $\|(\frac{1}{\alpha}
\phi^*\varphi^*)^k h\|_{sup} \leq (\frac{1}{\alpha})^k\|h\|_{sup}$,
it follows that the series given by the expression
$\sum_{k=0}^{\infty} ( \frac{1}{\alpha} \phi^*\varphi^*)^k h$,
converges absolutely to a bounded and continuous function $h^v$ on
$V(\bar{K}_v)$. Let $\|.\|_{\varphi,v}=\|.\|_1 \exp (h^v) $, then
$\|.\|_n$ converges uniformly to $\|.\|_{\varphi,v}$ and it is not
hard to check that $\|.\|_{\varphi,v}$ satisfies
$$\|.\|_{\varphi,v}=(\phi^*\varphi^*\|.\|_{\varphi,v})^{1/\alpha},$$
which was the result we wanted to prove.
\end{proof}
\begin{definition}
The metric $\|.\|_{\varphi,v}$ is called the canonical metric on
$\L_v$.
\end{definition}
\begin{example}
Consider the line bundle $\L=\cO_{\Pl^n}(1)$ on $\Pl_{\Qbar}^n$ and
the map $\phi_2(T_0:...:T_n)=(T_0^2:...:T_n^2)$. If we choose the
Fubini-Study metric $\|(\lambda_0 T_0+...+\lambda_n T_n
)(a_0:...:a_n)\|_{FS}=\frac{|\sum \lambda_i a_i |}{\sqrt{\sum_i
a^2_i}}$ as our smooth metric $\|.\|_{1}$ on $\L_{\C}$, the limit
metric we get at infinity is:
$$\|(\lambda_0 T_0+...+\lambda_n T_n)(a_0:...:a_n)\|_{nv}=\frac{|\sum
\lambda_i
a_i |}{\sup_i(|a_i|)}.$$
\end{example}

\subsection{Canonical measure and integral at infinite places}
In this part we set up what will be called the
integral at infinite places attached to a map $\varphi: V
\rightarrow V$ and a place $v$ of $K$ over infinity.  We need to
develop some analytic theory related to the complex points of $V$.\\
Let $M$ be a n-dimensional complex projective variety, $\vfi : M
\rightarrow M$ a map on $M$ and $\L$ an ample line bundle on $M$.
Suppose also that for some $\alpha
>1$ we have $\phi : \L^{\alpha} \xrightarrow{\sim} \varphi^*\L$, and
that the line bundle $\L$ is equipped with the canonical metric
$\|.\|_{\varphi}$ on the fibres.  Let $U \subset M$ be an open set.
The function $x \mapsto -\log \|s(P)\|_{\vfi}$ for a non-zero
holomorphic section $s$ on $U$, is not necessarily smooth, and due
to this fact, the first Chern ``form''
$c_1(\L,\|.\|)=\frac{1}{({\pi}i) }\partial\overline{\partial}\log
\|s_1(P)\|_{\vfi}$ may be no more than a distribution. We would like
to define the product $c_1(\L,\|.\|) \wedge...\wedge c_1(\L,\|.\|)$.
Unfortunately for general currents we do not have a product as we do
for smooth currents. Results of Bedford, Taylor and Demailly
\cite{bet}, \cite{dema1}, \cite{dema2} and \cite{dema3}, allow us to
consider a product of currents with good properties.  We follow the
presentation in \cite{maillot}.
\begin{definition} (Lelong).
Let $U$ be an open set of complex manifold $M$ of dimension $n$. A
current $T \in D^{p,p}(U)$ is said to be positive ($T \geq 0$) if
for every choice of $C^{\infty}$ $(1,0)$-forms
$\alpha_1,...,\alpha_{n-p}$ with compact support on $U$, the
distribution $T \wedge (i\alpha_1 \wedge \bar{\alpha}_1) \wedge
...\wedge (i\alpha_{n-p} \wedge \bar{\alpha}_{n-p})$ is a positive
measure on $U$.
\end{definition}

\begin{example}
A locally integrable function $u$ on $M$ is says to be
plurisubharmonic if the hessian $i \partial \bar{\partial} u = i
\sum \partial^2u/\partial z_j \partial z_m z_j \wedge \bar{z}_m \geq
0$ on $M$. For basic properties of plusibharmonic functions see for
example \cite{dema1} or \cite{dema2}.
\end{example}

\begin{example}
Let $Y$ be an algebraic $p-$cycle on $M$. The $(p,p)-$current
$\delta_Y$ of integration on $Y$, is a positive current on $M$.
\end{example}

\begin{definition}(Bedford-Taylor). Let $T$ be a positive closed current
of type $(p,p)$ and $u$ a plurisubharmonic function locally
bounded on $U$. We define the product $(dd^c u) \wedge
T=dd^c(uT)$.
\end{definition}
\begin{remark}
The product $Tu$ is well defined and in general we can define $(dd^c
u_1)\wedge (dd^c u_2)... (dd^c u_q)\wedge T=dd^c(u_1 \wedge (dd^c
u_2)...(dd^c u_q) \wedge T).$ \\ By prop 1.2 in \cite{dema3}, the
current $(dd^c u_1)\wedge (dd^c u_2)... (dd^c u_q) \wedge T$ is a
positive closed current of bidegree $(p+q,p+q)$.
\end{remark}
\begin{lemma}
Let $\L$ be an ample line bundle on $M$ and let $\vfi : M
\rightarrow M$ be a map with the property that for some $\alpha
>1$, there is an isomorphism $\phi : \L^{\alpha}
\xrightarrow{\sim} \varphi^*\L$.  Assume that $\|.\|_{\vfi}$ is the
canonical metric on $\L$. The function $P \mapsto -\log
\|s(P)\|_{\vfi}$ is plurisubharmonic on the open set $U=M - \Div(s)$
and therefore the current $i
\partial \bar{\partial} (-\log \|s(P)\|_{\vfi})$ is a positive
current on $U$.
\end{lemma}
\begin{proof}
The proof is basically taken from \cite{kawaguchi}. Consider the
continuous and positive function
$H(P)=\frac{\|s(P)\|_2}{\|s(P)\|_1}$ on $M$. Define $c=\min_{P \in
M} H(P)$. By changing $\vfi$ by $c \vfi$ if necessary, we can have
$H > 1$ (this will not affect the result by the second part of
theorem 2.2 in \cite{zhangadelic}). In general we have
$$\frac{\|s\|_{n}}{\|s\|_{n-1}}=\frac{\|s \circ \vfi \|_{n-1}}{\|s
\circ \vfi \|_{n-2}},$$ and
then $-\log \|s(P)\|_n \leq -\log
\|s(P)\|_{n-1}$ for every $n > 2$ and the sequence $\{-\log\|.\|_n
\}^{\infty}_{n=1}$ is a non increasing sequence of plurisubharmonic
function converging to $-\log \|s(P)\|_{\vfi}$.
\end{proof}

\begin{proposition}
Let $s_i$ $(i=1..q)$ be sections of the line bundles $\L_i$
respectively, such that the divisors $\Div(s_i)$ meet properly on
$M$. Let us denote $c_1(\L_i)=c_1(\L_i,\|.\|_{\vfi})$ the Chern
``form'' associated to the canonical metric studied in proposition
\ref{canonical metric}, then the current
$$(- \log
\|s_i\|_{\vfi})c_1(\L_1)...c_1(L_{i-1}).\delta_{\Div(s_{i+1})}...\delta_{\Div(s_{q})}
$$
is a well defined current and
$$\int_M (- \log
\|s_i\|_n)c_1(\L_1,\|.\|_n)...c_1(L_{i-1},\|.\|_n).\delta_{\Div(s_{i+1})}...\delta_{\Div(s_{q})}$$
tends to
$$\int_M (- \log
\|s_i\|_{\vfi})c_1(\L_1)...c_1(L_{i-1}).\delta_{\Div(s_{i+1})}...\delta_{\Div(s_{q})}.$$
\end{proposition}
\begin{proof}
We have that $c_1(\L_i,\|.\|_{\vfi})$ is a positive current that can
be written locally in the form $dd^c u$ where $u=-\log \|.\|_{\vfi}$
is a plurisubharmonic function on $M$.  On the other hand
$\delta_{\Div(s_i)}$ are closed and positive (Check theorem 3.5 in
\cite{dema1}).  The sequence of currents $$\{- \log
\|s_i\|_nc_1(\L_1,\|.\|_n)...c_1(L_{i-1},\|.\|_n).\delta_{\Div(s_{i+1})}...\delta_{\Div(s_{q})}
\}_n$$ converge weakly to $- \log
\|s_i\|_{\vfi}c_1(\L_1)...c_1(L_{i-1}).\delta_{\Div(s_{i+1})}...\delta_{\Div(s_{q})}$
after the following general proposition proved in \cite{dema3}.
\end{proof}
\begin{proposition} \label{demaillyconv}
(Demailly) Denote by $L(u)$ the set of point where the
plurisubharmonic function $u$ is not locally bounded. Let $U$ be an
open set of $M$ and $T \in D_+^{p,p}(U)$ a positive closed current
of type $(p,p)$. Let also $u_1,...,u_q$ be plurisubharmonic
functions on $U$, such that for every choice of indices
$j_1<j_2...<j_m$ inside $\{1,2,...,q\}$ the intersection $L(u_{j_1})
\cap...\cap L(u_{j_m}) \cap Supp(T)$ is contained in an analytic set
of complex dimension $\leq n-p+m$. One can then construct the
currents $u_1 (dd^cu_2) \wedge...\wedge (dd^c u_q) \wedge T$ and
$(dd^c u_1) \wedge...\wedge (dd^c u_q) \wedge T$ of mass locally
finite over $U$ and uniquely characterized by the fact: for every
non-increasing sequences $(u^k_1)$,...,$(u_q^k)$ of plurisubharmonic
functions converging punctually to $u_1,...u_q$ respectively, we
have that $u^k_1 (dd^cu^k_2) \wedge...\wedge (dd^c u^k_q) \wedge T$
and $(dd^c u^k_1) \wedge...\wedge (dd^c u^k_q) \wedge T$ converge
weakly on $U$ to $u_1 (dd^cu_2) \wedge...\wedge (dd^c u_q) \wedge T$
and $(dd^c u_1) \wedge...\wedge (dd^c u_q) \wedge T$ respectively.
\end{proposition}
\begin{proof}
For the proof we refer to \cite{dema3}, Thm. 3.4.5 and Pro. 3.4.9.
or \cite{dema2}. Thm. 2.5 and Pro. 2.9.
\end{proof}
\begin{definition} The canonical current
associated to $\vfi$ is defined as $T_{\vfi}=c_1(\L,\|.\|_{\vfi})$.
The canonical distribution associated to $\vfi$ is
$$d\mu_{\vfi}=c^n_1(\L_v,\|.\|_{\vfi}).$$
\end{definition}
\begin{proposition}
The canonical distribution is in fact a measure, which we call the
canonical measure.
\end{proposition}
\begin{proof}

The $(1,1)-$current $T_{\vfi}$ can be identified with an expression
$T_{\varphi}=\sum_{i,j}T_{i,j}dz_i \wedge d\bar{z}_j$ where the
coefficients $T_{i,j}$ are distributions.  Consider the
$(n-1,0)-$form $\alpha=\sum_{|I|=n-1} \alpha_I dz_{C_I}$.  The fact
that $T_{\varphi} \geq 0$ forces $ \sum T_{i,j} \alpha_{n-i}
\bar{\alpha}_{n-j}$ to be a positive measure for every $\alpha$.  As
a consequence the $T_{i,j}$ are complex measures with
$T_{i,j}=\bar{T}_{i,j}$.  In the same way the $p-$current
$\bigwedge^p_{i=1} T_{\vfi}$ has measures as coefficients for each
$p$.
\end{proof}

\begin{proposition}
Suppose that $A$ is a subset of $M$ such that $\mu(A)=\int_A
d\mu_{\varphi} < \infty$, then
\begin{enumerate}
\item $\mu_{
\varphi}(\varphi(A)) = (\deg \varphi)^n \mu_{\varphi}(A) < \infty$
whenever $\varphi|A$ is injective,
\item $\mu_{\varphi}(\varphi^{-1}(A))= \mu_{
\varphi} (A) < \infty .$
\end{enumerate}
\end{proposition}
\begin{proof}
Let $deg(\vfi)=d$ be the algebraic degree of $\vfi$.  Take an open set $W$
with $\mu(W) < \infty$. Assume that we have
$\varphi^{-1}(W)=\bigcup^{d^n}_{i=1} U_i$ where $\varphi : U_i
\rightarrow W$ is injective for each $i$ and let $U$ denotes any of
the $U_i$.  Consider $n$ local sections $s_i \neq 0$ of $\matO(1)$
holomorphic on $W$, in this case $d\mu_{\varphi}=\frac{1}{({\pi} i )
}\partial\overline{\partial}(g_1)...\frac{1}{({\pi} i )
}\partial\overline{\partial}(g_n)$, where
$g_i=\log(\|s_i(P)\|_{\varphi})$ and we have
\begin{equation*}
\begin{split}
\mu_{\varphi}(W) & = \int_{W}\frac{1}{({\pi}i)
}\partial\overline{\partial}\log(\|s_1(P)\|_{\varphi})...
\frac{1}{({\pi}i)
}\partial\overline{\partial}\log(\|s_n(P)\|_{\varphi}) \\
& = \int_{U} \frac{1}{({\pi}i)
}\partial\overline{\partial}\log(\|s_1(\varphi(P))\|_{\varphi})...
\frac{1}{({\pi}i)}\partial\overline{\partial} \log(
\|s_n(\varphi(P))\|_{\varphi})\\
& = d^n \int_{U} \frac{1}{({\pi}i)
}\partial\overline{\partial}\log(\|(\varphi^*s_1)(P)\|_{\varphi})...
\frac{1}{({\pi}i)}\partial\overline{\partial}\log(\|(\varphi^*s_n)(P)\|_{\varphi})\\
& = d^n \mu_{ \varphi}(U),
\end{split}
\end{equation*}
and
$$\mu_{\varphi}(\varphi^{-1}(W)) = d^n \mu_{ \varphi}(U)
= \mu_{ \varphi} (W),$$
which is the result we wanted to prove.
\end{proof}

\begin{definition}
Let $V$ be a projective variety defined over a number $K$, $\vfi : V
\rightarrow V$ a map on $V$ and  $(\L,\|.\|_v)$ a metrized ample
line bundle on $V$ with the property that there exist an isomorphism
$\psi : \L^{\alpha} \xrightarrow{\sim} \varphi^*\L$.  Let $v : K
\hookrightarrow \C$ be a place of $K$ over infinity.  The canonical
measure $d\mu_{\vfi, v}$ is the canonical measure on $V \otimes_v
\C$ associated to $\vfi_v=\vfi \otimes_v \C$ and $\L_v=\L \otimes_v
\C$.
\end{definition}

\subsection{Examples}
Here we revisit some examples of self maps on algebraic varieties
with extra information about canonical heights and measure.
\begin{example}
Suppose that we are working with a number field $K$ and
$\varphi=(p_0:...:p_n) : \Pl_K^n \rightarrow \Pl_K^n$ is a map on
the $n-$dimensional projective space over $K$. In general it is hard
to get a closed form for the iterate of such a map. In case that for
some natural $k$ we have $p_i(T_0,...,T_n)=T^k_i$ for each $0 \leq i
\leq n$, we obtain the so-called naive height on $\Pl_K^n$:
$$h_{nv}([t_0 : ... :t_n]) = \frac{1}{[K:\mathbb{Q}]}\log\prod _{\text{places
} v
  \text{ of } K} \sup(|t_0|_v ,...,|t_n|_v)^{\N_v},$$
where $\N_v = [K_v:\bQ_w]$ and $w$ is the place of $\bQ$ such that
$v \mid w$.  The associated measure $d\mu_{\vfi}$ is the normalized
Haar measure on the Torus $S^1 \times...\times S^1$.  If $T_0,
T_1,...,T_n$ represent projective coordinates in $\Pl^n$, the
canonical metric at infinity whose curvature gives the canonical
measure is
$$\| (\lambda_0 T_0 + ... + \lambda_n T_n) (a_0:...:a_n)\|_{nv} =
\frac{|\lambda_0a_0+...+\lambda_na_n|}{\sup (|a_0|,...,|a_n|)}.$$
\end{example}

\begin{example}
Let $A$ be an abelian variety and $\L$ a symmetric line bundle on
$A$. The multiplication by n, $[n] : A \rightarrow A$, satisfies
$[n]^* \L \backsimeq \L^{n^2}$, the canonical height is the
Neron-Tate height $\hat{h}_{NT}$ and the canonical measure, after
results in \cite{moretasterisque}, is the Haar measure on $A$.
\end{example}

\begin{example}
For the smooth toric variety $\Pl(\Delta)$ and the morphism $[p]:
\Pl(\Delta) \rightarrow \Pl(\Delta)$ we have $[p]^*
(\cO_{\Pl(\Delta)}) \backsimeq \cO_{\Pl(\Delta)}^p$ and we can build
a canonical height $\hat{h}_{\cO(\Pl(\Delta)),p}$.  In
\cite{maillot}, V. Maillot presents explicit formulaes for the
canonical measure as well as several properties of the canonical
height. The Mahler formula is established in this case as a
consequence of the vanishing of the canonical multiheight of the
whole variety $\Pl(\Delta)$.
\end{example}

\section{Arithmetic intersection theory}
In this section we develop the arithmetic intersection theory for
arithmetic varieties with Cohen-Macaulay local rings.  The
n+1-dimensional variety $X$ will be Macaulay and equipped with a
finite type, flat and proper map $f: X \rightarrow \Spec(\cO_K)$,
where $K$ is a number field.  First we introduce the geometric
intersection:
\begin{definition}
We say that a $q$-cycle $D$ in $X$ is locally regular complete
intersection (l.r.c.i.) if it is locally given by the intersection
of a regular sequence of length $n+1-q$.
\end{definition}
\begin{remark}
A l.r.c.i. $n$-cycle $D$ is just a Cartier divisor.
\end{remark}
\begin{definition}\label{schematic intersection}Suppose that the $q_1-$cycle
$D_1$ is
l.r.c.i. and $D_2$ is any $q_2-$cycle. Assume that they have no
common components.  Then the intersection is given by
  $$(D_1.D_2) =\sum^{q_1}_{i=0} (-1)^i
\Tor_i(\cO(D_1),\cO(D_2)).$$ In case $q_1+q_2 \leq n+1$ we define
the degree of the intersection as
$$\deg(D_1.D_2) =\sum^{q_1}_{i=0} (-1)^i \length
(\Tor_i(\cO(D_1),\cO(D_2)).$$
\end{definition}
\begin{definition}
If a cycle $D_1$ is such that for some natural m, $mD$ is l.r.c.i.,
we can extend the intersection as $(D_1.D_2)=1/m(mD_1.D_2)$.  An
n-cycle with this property is called a $\Q-$Cartier divisor.
\end{definition}
The intersection just defined satisfy many of the desirable
properties for intersection numbers.  Symmetry is clear from the
definition and associativity is a consequence of the convergence of
the $\Tor$ spectral sequence \cite{Serre}. Bilinearity is a
consequence of the lemma:
\begin{lemma}
Let $A$ be a commutative ring, $I$ an ideal and $f$ a non-zero
divisor in $A$.  Then we have the exact sequence:
$$0 \rightarrow A/I \rightarrow A/fI \rightarrow A/fA \rightarrow 0.$$
\end{lemma}
\begin{proof}
The canonical map from $A/fI \rightarrow A/fA$ has kernel $fA/fI$.
On the other hand the map $\gamma : A \rightarrow fA/fI$, given by
$\gamma(1)=[f]$ has kernel $I$, because if two elements $a \in A$
and $i \in I$ satisfy $af=if$, then $a=i$ because $f$ is not a
zero divisor.
\end{proof}
\begin{proposition}\label{linear equivalence}
     Let $X$ be a projective arithmetic variety, $C$ a Cohen-Macaulay
    projective curve in $X$ and $D$ a Cartier divisor, then $$
    \deg(D.C) = \deg_C\mathcal{O}_X(D)_{|C}.$$
  \end{proposition}
  \begin{proof}
   If $C$ is a projective curve and $L$ is a line bundle on $X$ a
  projective variety, one can then speak of $(L.C)$ for $L$ is the
  difference in $\Pic(X)$ between two very ample line bundles each of
  them having sections with no common components with $C$.  The
  result follows because a line bundle on a Cohen Macaulay curve has
  a well defined degree.
\end{proof}
\begin{definition}
 If
$\L$ is a line bundle on $X$ we will denote by $c_1(\L)$ the class
of divisors determine by $\L$ and by $c_1(\L)^i$ the intersection of
any element in this class with itself i times.
\end{definition}

\begin{proposition}\label{codimension}
    If $D_2$
has codimension $n+1$ and $D_1$ is locally given by one equation,
then $\deg(D_1.D_2) = 0$.
  \end{proposition}
  \begin{proof}
  Suppose that $f$ is the equation defining $D_1$ and $I$ is the
  ideal of $D_2$ in the local ring $\cO_x=A$. The assumption on the
dimension of $D_2$ implies
  that the modules $A/I$ and $\Tor_i(A/I,A/f)$ are of finite length for all
i. Now the
  result follows because the length is an additive function and we have the
exact sequence:
$$0 \rightarrow \Tor_1(A/I,A/f) \rightarrow A/I
  \rightarrow A/I  \rightarrow A/(I+(f)) \rightarrow 0.$$

\end{proof}
\begin{proposition} \label{contraction} Suppose that $\sigma : X_1
\rightarrow X$
is a map of projective Arithmetic varieties over $\Spec(\cO_K)$. Let
$C$ be a closed $1-$cycle of $X_1$ and $\L$ a line bundle on $X$. If
$C$ is
  contracted by $\sigma$ to a subscheme of $X$ of codimension n+1 the
  intersection number $\deg(\sigma^*(L).C)$ is zero.
  \end{proposition}
  \begin{proof}
$L$ can be realized as the line bundle associated to
  the difference of two very ample divisors on $X$ each of them having
  no intersection with $\sigma_*(C)$. The reciprocal images of these
  divisors in $Y$ do not meet $C$ and the result follows.
  \end{proof}

The geometric intersection, however,
does not take into account the places of $K$ over infinity. An
original ideal of Arakelov \cite{arakelovinter}, later developed by
Szpiro \cite{asterisquedeux}, Bost \cite{BGS}, Gillet \cite{BGS},
Soul\'{e} \cite{BGS}, Zhang \cite{zhangthese} \cite{zhangvar},
Faltings and many others, allow us to consider an intersection
theory that equally value all places of $K$.
\begin{definition}
Let $X$ be a Cohen-Macaulay arithmetic variety of dimension $n+1$,
defined over a number field $K$ and $\bar{\L}=(\L,\|.\|)$ a
hermitian line bundle on $X$.  Suppose that $v_1,...,v_s$ are the
different places of $K$ over infinity.  We will denote by
$\hat{c}_1(\L)$ the vector
$(c_1(\L),c_1(\bar{\L})_{v_1},...,c_1(\bar{\L})_{v_s})$, where
$c_1(\L)$ is the class of equivalent divisors determined by $\L$ and
$c_1(\bar{\L})_{v_i}$ is the $(1,1)-$current associated to the
metric $\|.\|_{v_i}$ on $\L_{v_i}$.
\end{definition}
\begin{proposition} \label{arithmetic intersection}Let $X$ be a
Cohen-Macaulay arithmetic variety of dimension $n+1$, defined
over a number field $K$.  Let $Z \in Z_k(X)$ be a cycle on $X$ and
$\bar{\L}_1,...,\bar{\L}_k$ a set of hermitian line bundles on $X$.
Then the number $\hat{\deg}_Z (\hat{c}_1(\L_1)...\hat{c}_1(\L_k))
\in \bR$, is completely determined by the properties:
\begin{enumerate}
\item{is $k-$linear.}
\item{is symmetric.}
\item for $k=0$ and $Z=\sum_i n_iP_i$ ($P_i \in X_v$), we have then
$\hat{\deg}_Z=\sum_i n_i
N_{P_i} \log N(v)$ where $N_{P_i}=[K(P_i):K]$.
\item for $k \geq 1$ and $s_k \neq 0$ a section of $\L_k$ which meets
$Z$ properly we have
\begin{equation*}
\begin{split}
\hat{\deg}_Z(
\hat{c}_1(\L_1)...\hat{c}_1(\L_{k})|Z)& = \hat{\deg}_Z(
\hat{c}_1(\L_1)...\hat{c}_1(\L_{k-1})|Z.\Div(s_k))\\ & -
\sum_{v/\infty} \int_{X(\bC)}
\delta_{Z(\C)} \log \|s_k\|_{k,v}
c_1(\bar{\L}_1)_v...c_1(\bar{\L}_{k-1})_v,
\end{split}
\end{equation*}
where $\sum_{v/\infty}$ is taken over the places of $K$ at infinity.
\end{enumerate}
\end{proposition}
\begin{proof}
Conditions (i), (iii) and (iv) are sufficient to determine
recursively the number
$\hat{\deg}_Z(\hat{c}_1(\L_1)...\hat{c}_1(\L_k))$. Suppose that we
consider sections $s_i$ of the line bundles $\L_i$ respectively,
such that the divisors $\Div(s_i)$ meet properly in $X$. Let $v$ be
a place of $K$ over infinity.  Introducing the star product
$g_{1,v}*g_{2,v}*...*g_{k,v}$ of the currents $g_i=-\log\|s_i\|_v$
we can state a non-recursive formula of the arithmetic degree
$$\hat{\deg}_Z(\hat{c}_1(\L_1)...\hat{c}_1(\L_k))=
(\Div(s_1)...
\Div(s_k))_{\finite}+\sum_{v/\infty}\int_{Z(\C)}g_{1,v}*...*g_{k,v},$$
where the first term on the right is representing the weighted sum
$$(\Div(s_1)....
\Div(s_k))_{\finite}=\sum_{v \finite} \deg(\Div(s_1)_v...
\Div(s_k)_v) \log N(v).$$

The condition (ii) will be a consequence of the following
lemma.
\end{proof}
\begin{lemma}The arithmetic degree $\hat{\deg}_Z
(\hat{c}_1(\L_1)...\hat{c}_1(\L_k))$ is a symmetric function of the
$\L_i$.
\end{lemma}
\begin{proof}
The geometric intersection is symmetric on the divisors $\Div(s_i)$.
Let's concentrate then in the term involving Green functions over a
fixed place $v$.  Suppose that $g_i$ (for $i=1,2$) are Green
currents of "log" type along the cycles $Z_1$ and $Z_2$ and relative
to $v$ (for the existence see lemma 1.2.2 of \cite{BGS}). We have
$g_1*g_2=g_2*g_1+\delta T_1+ \bar{\delta} T_2$ for some currents
$T_1 \in D^{p-1,p}$ and $T_2 \in D^{p,p-1}$, then $\int_Z g_1*g_2=
\int_Z g_2*g_1 + \int_Z \delta T_1 + \int_Z \bar{\delta} T_2$ and by
Stokes theorem we obtain the symmetry for the arithmetic degree.
\end{proof}
\begin{proposition}
Suppose that $\L$ is a hermitian line bundle on $X$ and $f \in K(X)$
is a rational function on $X$. Then: $$\hat{\deg} (\hat{c}_1
(\cO(f)). \hat{c}^n_1(\L))=0.$$
\end{proposition}
\begin{proof}
The curvature of the trivial bundle $c_1(\cO(f))=0$. Using this
result and the symmetry of the arithmetic degree we can reduce to
the case of dimension 1, which is nothing else but the product
formula (see for example the treatment in \cite{orsay}).
\end{proof}

\begin{proposition} \label{inters-at-infinity}
Suppose that $\L$ is a hermitian line bundle on $X$ and $f \in K(X)$
is a rational function on $X$ then,
$$\hat{\deg} ( \hat{c}^n_1(\L)|\Div(f))= \sum_{v/\infty} \int_{X_{v}}
\log|f|_{v}
d\mu_{v}.$$
\end{proposition}
\begin{proof} We have that $\hat{\deg} (\hat{c}_1 (\cO(f)).
\hat{c}^n_1(\L))=0$ and also that
\begin{equation*}
\begin{split}
\hat{\deg} (\hat{c}_1 (\cO(f)). \hat{c}^n_1(\L)) = \hat{\deg}(
\hat{c}_1(\L)^{n})|\Div(f))- \sum_{v/\infty} \int_{X_{v}(\bC)} \log
|f|_{v} (c_1(\bar{\L})^n.
\end{split}
\end{equation*}
which gives the formula we wanted.
\end{proof}
\begin{remark}In the notation of
\cite{sz-t-p} we can write $$(div(f). \L)_{\Ar}=\hat{\deg} (
\hat{c}^n_1(\L)|\Div(f)).$$
\end{remark}
\begin{definition}
Let $Y \in Z_q(X)$ be a $q-$cycle inside the arithmetic variety $X$
and $\L_1,...\L_q$ ample line bundle on $X$. The real number
$$h_{\L_1,...,\L_q}(Y)= \hat{\deg}_Y
(\hat{c}_1({\L_1})...\hat{c}_1({\L_q})|Y)$$
is called the multi-height of $Y$ relative to $\L_1,...,\L_q$.
\end{definition}
Let's denote by $X_K$ the generic fibre of $X$.  Suppose that we
have a map $\vfi : X_{K} \rightarrow X_{K}$ and an ample line bundle
$\L$ on $X_{K}$ such that for some number $\alpha$, we have the
isomorphism of line bundles $\psi : \L^{\alpha} \cong \vfi^* \L $.
Then we can build a sequence of arithmetic varieties $X_k$,
$k=1,2,...$, models of $X_{K}$ over $\Spec(\cO_K)$, and metrized
line bundles $(\L_k,\|.\|_k)$ on $X_k$, such that
$\|.\|_{k+1,v}=(\psi^* \vfi^* \|.\|_{k,v})^{1/\alpha}$ for each
place $v$.  The detailed discussion of the construction of the
$\L_k$ can be found in page 10 and 11 of \cite{zhangadelic}, it
constitutes an example of a sequence of adelic metrized line
bundles, that is, a sequence of line bundles with good metrics at
every place.  The numbers $\hat{\deg}_Y
(\hat{c}_1((\L_{k},\|.\|_k))^{\dim(Y)+1}|Y)$ converge for every
$p-$cycle $Y \subset X_K$ (see theorem 1.4 of \cite{zhangadelic})
and the limit will be called
$\hat{\deg}({c_1(\L,\|.\|_{\vfi})}^{\dim Y+1}|Y)$.  This number is
now not depending on the $X_k$.
\begin{definition}  \label{ canonical multiheight} Under the conditions just
discussed we define
$$h_{\vfi}(Y)=\frac{\hat{\deg}({c_1(\L,\|.\|_{\vfi})}^{\dim Y +
1}|Y)}{c_1(\L)^{\dim Y}}.$$
\end{definition}
\begin{proposition}
Consider a rational function $F$ on $X_K$ and $\L$ a line bundle
on $X_K$ with conditions as before.  Then, If the map $\vfi :
X_{K} \rightarrow X_K$ extends to a map $\vfi : X \rightarrow X$,
we have
$$h_{\vfi}(\Div(F))= \sum_{v/\infty} \int_{X_v} \log |F|_v
\frac{d\mu_{\vfi,v}}{c_1(\L)^n}.$$
\end{proposition}
\begin{proof}
This is a consequence of proposition \ref{inters-at-infinity} and
definition \ref{ canonical multiheight}.
\end{proof}
\section{The blow up}
Let $\vfi: \Pl_K^n \too \Pl_K^n$ be a map with $\vfi^* \cO(1) \cong
\cO(d)$ and defined over a number field $K$. A model for $\vfi$ over
$\cO_K$ is a vector $(q_0,q_1,...,q_n)$ of elements of
$\cO_K[T_0,T_1,...,T_n]_d$, such that our map is expressed as
$\vfi=(q_0:q_1:...:q_n) : \Pl_K^n \too \Pl_K^n$. We will be
interested in models of $\vfi$ such that $(q_0,q_1,...,q_n)$ form a
regular sequence.
\begin{lemma}
Let $A=\cO_K[T_0,T_1,...,T_n]$ and $p_0,p_1,...,p_n \in A$. Then,
the following two statements are equivalent:
\begin{enumerate}
\item{the sequence $p_0,p_1,...,p_n$ is a regular sequence.}
\item{$\dim(V( \langle p_0,p_1,...,p_n \rangle))=0$.}
\end{enumerate}
\end{lemma}
\begin{proof}
The ring $A$ is Cohen-Macaulay of dimension $n+1$, therefore the
sequence $p_{0},...,p_{n}$ is regular if and only if it is a maximal
system of parameters.
\end{proof}
\begin{corollary}
Let $A=\cO_K[T_0,T_1,...,T_n]$ and $p_0,p_1,...,p_n$ a regular
sequence in $A$. Then the sequence
$p_k=\{(p_{k0}:p_{k1}:...:p_{kn})\}_k$, defined recursively by
$$p_0=(p_{0}:p_{1}:...:p_{n}) \quad p_{ki}=p_{k-1i}(p_{0}:p_{1}:...:p_{n})
\quad 0 \leq i \leq n \quad  k > 0$$ is also a regular sequence
for all $k$.
\end{corollary}
\begin{proof}
If $p_{k-1}$ is a regular sequence, $\dim(V( \langle
p_{k-1,0},p_{k-1,1},...,p_{k-1,n} \rangle))=0$, because $\varphi$ is
a finite map, $\dim(V( \langle p_{k,0},p_{k,1},...,p_{k,n}
\rangle))=0$ and $p_k$ is also a regular sequence.
\end{proof}

\begin{definition} Denote $\Pl^n_{\cO_K}$ by $X$ and by $Y_k$, the closed
subscheme
of $X$ defined by the ideal $I_{k}=\langle p_{k,0},p_{k,1},...,p_{k,n}
\rangle$.  The model $X_k$ is defined by the property that
$\sigma^k:X_k$ $\rightarrow$ $X=\Pl^n_{\cO_K}$ is the blowing up
of $Y_k$. The exceptional divisor will be denoted by $E_k$ and its
irreducible components by $C_{v,i,k}$, in such a way that we have
a finite sum $E_k=\sum_{v,i>0}r_{v,i,k} C_{v,i,k}$.
\end{definition}

In the rest of this subsection we will work with $X=\Pl^n_{\cO_K}$ and a map
$\vfi : \Pl_K^n \rightarrow \Pl_K^n$
represented by a regular sequence $(p_0,...,p_n)$ in $\cO_K$.  In the way we
have defined the map $\sigma^k$, we have
$\sigma^kI_k = \cO_{X_k}(-E_k)$, by the universal property of the blow-up.
By the same property, the surjection
$\cO_{X_k}^2 \twoheadrightarrow \sigma^{k*}(\cO_X(d)) \otimes
\cO_{X_k}(-E_k)$, gives rise to a map $\vfi_k : X_k \rightarrow X$.
By definition of the map $\varphi_k$ we have $\vfi_k^*
\cO_X(1)=\sigma^{k*}(\cO_X(d)) \otimes \cO_{X_k}(-E_k)$.  We will denote
by $\L_0$ the line bundle $\cO(1)$ on $X=\Pl^n_{\cO_K}$ and $\L_k=\vfi^{k*}
\L_0$ on the model $X_k$.

\begin{proposition} \label{general blow-up}
The scheme $X_k$ is Macaulay and $Y_k$ is a subscheme of codimension
$n+1$ in $\cO_K$ and does not meet the generic fiber $X_K$. Each
component $C_{v,i,k}$ is isomorphic to the projective space of
dimension n over the the residual field $K_{v,i,k}$ of the close
point image of $C_{v,i,k}$.  If $X_{k,v}$ does not meet $Y_k$, then
$X_{k,v}$ is isomorphic to $\Pl^n_{K_v}$.
\end{proposition}
\begin{proof}
Since $I_k$ is finitely generated the scheme $X_k=\proj(\bigoplus_{n \geq 0}
I_{k}^n)$ is locally complete intersection in $\Pl^n_{\cO_K}$ and therefore
Macaulay.  The sequence $(p_{k,0},p_{k,1},...,p_{k,n})$ being regular in
$\cO_K[T_0,T_1,...,T_n]$ forces $Y_k$
to be of codimension $n+1$. By
definition of the blow-up we get that the components $C_{v,i,k}$ are
isomorphic respectively to $\Pl^n_{K_{v,i,k}}$.  The total fiber
$F_v$ and $E_k$ are Cartier divisors, therefore the components
$C_{v,i,k}$ are $Q-$Cartier divisors, because they don't meet each
other.  The intersection of $\L_k$ with itself n times
${c}_1(\L_k)^n$ represent a Cohen Macaulay curve on $X$. Using
propositions \ref{contraction} and \ref{linear equivalence} we have
\begin{equation*}
\begin{split}
\hat{\deg}({c}_1(\L_k)^n|C_{i,v,k}) &
=\hat{\deg}({c}_1(\cO_{X_k}(-E_k))^n|C_{i,v,k})
\\ &={\deg}(\cO_{\Pl^1_{K_{v,i,k}}}(1))=[K_{v,i,k}:K_v] \log|N(v)|.
\end{split}
\end{equation*}
Also for places $v$ of good reduction,
$$\hat{\deg}({c}_1(\L_k)^n|X_{k,v})=\deg(\cO_{\Pl^1_{K_v}}(1))=
\log|N(v)|.$$ In this way the proposition gives us a way to compute
the arithmetic intersection of $\L_k$ with the different vertical
components of $X_k$.
\end{proof}

\begin{definition}
Recall that $f : X \rightarrow \Spec(\cO_K)$ is an Arithmetic
variety.  The projection $f (Y_1) \subset \Spec(\cO_K)$ will be
called the places of bad reduction of $\varphi$.
\end{definition}
\begin{remark}
The only places $v$ appearing in the exceptional divisors
$E_k=\sum_{v, i>0}r_{v,i,k} C_{v,i,k}$ are the places of bad
reduction.
\end{remark}

Let $F_0 \in \cO_K[T_0,...,T_n]$.  We can assume that $v(F_0)=0$ for
every place $v$ of bad reduction because there is only finitely many
places of bad reduction and any Dedekind domain with finitely many
primes (like $\cap_{v} \cO_v(z)$) is unique factorization domain.
Now consider the rational function $F=F_0/T^{\deg(F_0)}_n$ on
$\Pl^n_{\cO_K}$ and $F_k=\sigma^*_k F$.  The symbol $\infty_k$ will
be denoting the divisor of $X_k$ defined by the equation
$\sigma_k^*T_n=0$ and in particular $\infty=\Div(T_n)$.  If we
define the irreducible horizontal divisor $D$ in $X$ by the equation
$\Div(F) = D - \deg(F) \infty + \sum_{\finite \,v} v(F) X_v$, we can
establish the following lemma.

\begin{lemma} \label{Fpullback}
There exist non-negative integers
$x_{v,i,k}$ and $y_{v,i,k}$ depending only on $D$, such that $\Div
(F_k)$ can be written as
\begin{equation} \label{Mahler-Step-k}
\begin{split}
\Div(F_k) = D_k & - \deg(F) \infty_k + \sum_{v,i} x_{v,i,k}
C_{v,i,k} \\ & -\deg(F) \sum_{v,i} y_{v,i,k} C_{v,i,k} +
\sum_{\finite \,v} v(F) X_{v,k}
\end{split}
\end{equation}
where $D_k$ is the proper transform of $D$ by $\sigma_k$.
\end{lemma}
\begin{proof}
We have the formula for the divisor $\Div(F_k)$: $$\Div F_k =\Div
\sigma^*_k F = \sigma^*_k \Div F = \sigma_k^*(D)- n \sigma_k^*(
\infty) + \sum_{\finite \,v} v(F) \sigma_k^*(X_{v,k}),$$ but now,
for certain non-negative integers $x_{v,i,k}$ the reciprocal image
of the effective divisor $D$ is $\sigma_k^*(D)=D_k + \sum_{v,i}
x_{v,i,k} C_{v,i,k}$. Now also for certain non-negative integers
$y_{v,i,k}$ we have $\sigma_k^*(\infty)= \infty_k + \sum_{v,i}
y_{v,i,k} C_{v,i,k}$ and the proof is finished.
\end{proof}
\begin{corollary}
With the notation as before we have $$\sigma_k^*(\infty)= \infty_k +
\sum_{v,i} y_{v,i,k} C_{v,i,k} \quad \sigma_k^*(D)=D_k + \sum_{v,i}
x_{v,i,k} C_{v,i,k}.$$
\end{corollary}
\subsection{Negativity conditions} It is interesting to look at a particular
type of
models.
\begin{definition}
We say that a model $(p_0,...,p_n)$ satisfies negativity conditions
if we have $\rad ( \langle p_{k,0},p_{k,1},...,p_{k,n}, T_n\rangle
)=\langle T_0,T_1,...,T_n \rangle$ in $\cO_K[T_0,...,T_n]$, for
every $k$.
\end{definition}
\begin{example}
This condition is satisfied for example if we are considering a map
$\varphi=(p_0:p_1) : \Pl^1 \rightarrow \Pl^1$, on the Riemann sphere
and $p_0$ is a monic polynomial in the variable $T_0$.
\end{example}
\begin{lemma}
If the model have negativity conditions the proper transform $\infty_k$ of
$\infty$ in $X_k$ is equal to the reciprocal image
$\sigma^k(\infty)$.
\end{lemma}
\begin{proof} \label{do not meet infinity}
It is enough to show that $\infty$ does not meet $Y_k$ and this is a
consequence of the fact that the ideal $\rad ( \langle
p_{k,0},p_{k,1},...,p_{k,n}, T_n\rangle )=\langle T_0,T_1,...,T_n
\rangle$.
\end{proof}

\begin{lemma}
If the model has negativity conditions, there exist non-negative
integers $x_{v,i,k}$ depending only on $D$, such that $\Div (F_k)$
can be written as
$$\Div(F_k) = D_k - \deg(F) \infty_k + \sum_{v,i,k} x_{v,i,k} C_{v,i,k} +
\sum_{\finite \,v} v(F) X_{v,k}$$
where $D_k$ is the proper transform of $D$ by $\sigma_k$.
\end{lemma}
\begin{proof}
The results follows by Lemma \ref{do not meet infinity} and Lemma
\ref{Fpullback}.
\end{proof}
\section{Finite places}
Let $v$ denote a valuation on $K$, $\cO_v$ the set of elements $z
\in K$ such that $v(z) \geq 0$, $x$ will denote a vector
$x=(x_0,...,x_n) \in \bar{K}^n$ and $P \in \Pl^n(\bar{K})$ a point
in the $n$-projective space over $\bar{K}$. The valuation $v$ is
assumed to be extended to an algebraic closure $\bar{K}$ of $K$. For
a vector $x=(x_0,x_2,...,x_n)$ we define $v(x)=min_i \{ v(x_i) \}$.
For a polynomial $p \in \bar{K}[T_0,...,T_n]$ we take $v(p)$ as the
valuation of the vector formed by its coefficients. For a sequence
of polynomials $(p_0,...,p_n)$ we put $v(p_0,...,p_n)=\min_i \{
v(p_i) \}$. Suppose that $\varphi : \Pl_K^n \rightarrow \Pl_K^n$ is
a rational map of algebraic degree $d$ given by homogeneous
polynomials $(p_0:...:p_n)$ over $\cO_v$, then we can define a map
$$S_v : \bar{K}^{n+1}-(0,...,0) \times \cO_v[T_0,...,T_n]_d^{n+1}
\rightarrow \R_{\geq 0}$$
$$S_v(x,(p_0,...,p_n))=v(p_0(x),...,p_n(x))-v(x_0^d,...,x_n^d)$$

The map $S_v,$ is in fact a well defined map $S_v : \Pl^n \times
\cO_v[T_0,...,T_n]_d^{n+1} \rightarrow \R_{\geq 0}$, which we still
denote by $S_v(P,(p_0,...,p_n))$. To see this, take any two sets of
homogenous coordinates for $P$, say $(x_0,...,x_n)$ and
$(y_0,...,y_n)=\lambda(x_0,...,x_n)$, then the valuation
$v(p_0(\lambda x),...,p_n(\lambda x))=dv(\lambda) +
v(p_0(x),...,p_n(x))$ and $v(\lambda^d x^d)=d v(\lambda) + v(x^d)$,
and the result follows.

\begin{definition}Suppose that the polynomial $F$ has divisor
$\Div(F)=D- \deg(F) \infty - \sum_{\finite v} v(F) X_v$, then we
define:
\begin{equation*}
\begin{split}
E(F, v \finite)= & -\limsup_k \sum_v \log|N(v)| \Big( \frac{\sum_{P
\in D} S_v(P ,p_{k,0},...,p_{k,n})}{d^{nk}}
\\ & -\deg(F) \frac{\sum_{P \in \infty}
S_v(P ,p_{k,0},...,p_{k,n})}{d^{nk}} - v(F) \Big).
\end{split}
\end{equation*}
\end{definition}
\begin{remark}
We have $S_v(P ,p_{k,0},...,p_{k,n}) > 0$ only for finitely many
$P$, because the sequence $(p_{k,0},...,p_{k,n})$ is regular.  In
dimension one we can actually change the $\limsup$ of the formula
into a $\lim$.
\end{remark}
\subsection{Convergence of each v-adic integral in dimension
one} This part basically follows section 5.1 of \cite{sz-t-p}.
Suppose that we are working in dimension one, i.e. with a map
$\vfi=(p_0:p_1) : \Pl^1 \rightarrow \Pl^1$ on the Riemann Sphere.
Let's keep the notation of the previous subsection, that is
$S_v(x,(p_0,p_1))=v(p_0(x),p_1(x))-v(x_0^d,x_1^d)$.
\begin{proposition} \label{convergence of finite place int}
The sequence
$$h_k(P)=\{\frac{S_v(P,(p_{k,0},p_{k,1}))}{d^k}\}_k$$
\begin{enumerate}
\item is bounded and increasing, and therefore convergent to a
function, which we denote $h_{p_0,p_1,v}(P)$. \item
$h_{p_0,p_1,v}(\varphi(P))=dh_{p_0,p_1,v}(P).$
\end{enumerate}
\end{proposition}
The proof will be the result of the application of two lemmas:
\begin{lemma} \label{S equation}
Suppose that we denote $P_k=(p_{k,0}(P):p_{k,1}(P))$, then we have
the equality
$$S_v(P,(p_{k+1,0},p_{k+1,1})) = d
S_v(P,(p_{k,0},p_{k,1}))+ S_v(P_k,(p_0,p_1)).$$
\end{lemma}
\begin{proof}
Assume that $P=(x_0:x_1)$ and $v(x)=v((x_0,x_1))=0$. If we set
$x_k=(p_{k,0}(x),p_{k,1}(x))$ we have the equalities
\begin{equation*}
\begin{split}
S_v(P,(p_{k+1,0},p_{k+1,1})) & =v(p_{k+1,0}(x),p_{k+1,1}(x)) \\
&= v(p_{0}(x_k),p_{1}(x_k))-v(x^d_k)+ v(x^d_k)\\
&= d S_v(P,(p_{k,0},p_{k,1}))+ S_v((P_k,(p_0,p_1)),
\end{split}
\end{equation*}
which gives the result we were trying to prove.
\end{proof}
\begin{lemma} The function $S_v(P,(p_0,p_1))$ is bounded
on $\Pl^1(\bar{K})$, so we can define
$$R_v(p_0,p_1)=\sup_{P \in \Pl^1} \{S_v(P,(p_0,p_1)) \}.$$
\end{lemma}
\begin{proof}
There exist elements $0 \neq b_i \in \cO_v$, where $0 \leq i \leq
1$, such that $b_ix^{d}_i \equiv 0 (\langle p_0,p_1 \rangle)$. If
$P=(x_0:x_1) \in \Pl^1$, with $x_i \neq 0$, then $S_v(P,(p_0,p_1))
\leq v(b_i)$. So in general $\sup_{P} \{S_v(P,(p_0,p_1)) \} \leq
\sup_i \{v(b_i) \}$.
\end{proof}
Now we can proceed to prove proposition \ref{convergence of finite
place int}.
\begin{proof}
From lemma \ref{S equation}, we get $0 \leq h_{k+1}(P)-h_k(P) \leq
R_v(p_0,p_1)/d^{k+1}$, so $\{h_k(P)\}$ is bounded by
$R_v(p_0,p_1)/(d-1)$ and therefore converges. On the other hand
\begin{equation*}
\begin{split}
h_{k}(\varphi(P))-dh_k(P) & =\frac{S(P_1,(p_{k,0},p_{k,1}))}{d^k}
- \frac{dS(P,(p_{k,0},p_{k,1}))}{d^k}\\
&=\frac{S(P,(p_{k+1,0},p_{k+1,1}))-v(x_1)}{d^k}
- \frac{dS(P,(p_{k,0},p_{k,1}))}{d^k}\\
&=\frac{S(P_k,(p_{0},p_{1}))-v(x_1)}{d^k} \leq
\frac{R_v(p_{0},p_{1})-v(x_1)}{d^k}.
\end{split}
\end{equation*}
By passing to the limit we get
$h_{v,p_0,p_n}(\varphi(P))=dh_{v,p_0,p_1}(P)$.
\end{proof}
\begin{definition} \label{definition of int at finite place}
Suppose that $\vfi : \Pl_K^1 \rightarrow \Pl_K^1$. We will define
the ``integral'' of $\log|F|_v$ over the finite place $v$ of a
polynomial $F=K[z_1,...,z_n]$ as
\begin{equation*}
\begin{split}
\int_{\Pl^1_{\bar{K}_v^n}} \log|F|_v d\mu_{v,\varphi}= &
\log|N(v)|
\Big( \sum_{P \in D} h_{v,p_0,p_1}(P)
\\ & -\deg(F) \sum_{P \in \infty} h_{v,p_0,p_1}(P) - v(F) \Big) ,
\end{split}
\end{equation*}
and then
\begin{equation*}
\begin{split}
E(F, v \finite) = \sum_v \int_{\Pl^1_{\bar{K}_v^n}} \log|F|_v
d\mu_{v,\varphi} .
\end{split}
\end{equation*}
\end{definition}

\subsection{Geometry of $E(F, v \finite)$} We want to relate $E(F, v
\finite)$ with the geometry of the blow-up.
Suppose that $\sigma_k: (\Pl^n)_k \rightarrow \Pl^n$ is the blow-up
associated with the model $(p_{k,0}:...:p_{k,n})$. Writing
\begin{equation} \label{blow-up}
\begin{split}
\sigma^*_kD=D_k+\sum_{v,i} x_{v,i,k}C_{v,i,k} &
\qquad \sigma_k^*(\infty)= \infty_k + \sum_{v,i} y_{v,i,k}
C_{v,i,k}.
\end{split}
\end{equation}
where $C_{v,i,k}$ are the different components of the exceptional fibre
above $v$, and $K_{v,i,k}$ denotes the field of definition of the
close point corresponding to $C_{v,i,k}$, we can state:
\begin{proposition} For every v  finite
place of $K$, we have:
\begin{equation*}
\begin{split}
\sum_i x_{i,v,k}[K_{v,i,k}:K_v] & =\sum_{P \in D}
S_v(P,p_{k,0},...,p_{k,n}),\\
\sum_i y_{i,v,k}[K_{v,i,k}:K_v] & = \sum_{P \in \infty}
S_v(P,p_{k,0},...,p_{k,n}).
\end{split}
\end{equation*}
\end{proposition}
\begin{proof}
Let $\sigma_{\cO_v,k} : (X_k)_v \rightarrow X_{\cO_v}$ be the
blow-up $\sigma_k : X_k \rightarrow X$, composed with the base
extension to $\Spec(\cO_v)$.  Let $D_v$ be the localization of $D$
at $v$.  We have the equation $\sigma^*_{\cO_v,k}(D_v)=D_{v,k} +
f_{v,k} C_{v,k}$ for some horizontal divisor $D_{v,k}$ in $X_k$ and
non-negative integers $f_{v,k}$.  Suppose that $J_v$ is the ideal
sheaf of $D_v$ in $X_{\cO_v}$, then $\sigma_{\cO_v,k}(J_v)$ will
correspond to the ideal sheaf of $D_{v,k} + f_{v,k} C_{v,k}$.
We are going to assume that there is only one $P=(a_0:...:a_n) \in X_v \cap
D$, because otherwise we will blow up one point at a time.  The subscheme of
$(X_k)_v$ determined by
$\sigma_{\cO_v,k}(J_v))$ is isomorphic to
\begin{equation*}
\begin{split}
& Proj(\cO_v[T_0,...,T_n]/\langle
({p}_{k,i}(a_0,...,a_n)T_j-{p}_{k,j}(a_0,...,a_n)T_i)_{i,j} \rangle)
\\ \cong &  Proj(\cO_v[T_0,...,T_n]/ \pi^{r_{k,v}} \langle
(\bar{p}_{k,i}(a_0,...,a_n)T_j-\bar{p}_{k,j}(a_0,...,a_n)T_i)_{i,j}
\rangle)
\\ \cong  &  Proj((\cO_v/\pi^{r_{k,v}}\cO_v)[T_0,...,T_n]) \cup \Spec(\cO_v)
\\ = & r_{k,v,P} \Pl^n_{K_{v}} \cup \Spec(\cO_v) .
\end{split}
\end{equation*}
where
$\bar{p}_{k,i}(a_0,...,a_n)={p}_{k,i}(a_0,...,a_n)/\pi^{r_{k,v,P}}$
for all $0 \leq i \leq n$ and the valuation
$r_{k,v,P}=S_v(P,p_{k,0},...,p_{k,n})$. \\
On the other hand $\sigma_{v,k}^* P_{v,i,k} \cong x_{v,i,k}
\Pl^n_{K_{v,i,k}} \cong x_{v,i,k} [K_{v,i,k}:K_v] \Pl^n_{K_{v}}$. As
a consequence of this two facts
$x_{i,v,k}[K_{v,i,k}:K_v]=S_v(P,p_{k,0},...,p_{k,n})$ and
$$\sum_i x_{i,v,k}[K_{v,i,k}:K_v]=\sum_{P \in D}
S_v(P,p_{k,0},...,p_{k,n}).$$ The second part is analogous using
$\deg(F)\infty$ instead of $D$.
\end{proof}
\begin{remark}
The expression $E(F, v \finite)$ previously defined, takes the
form
\begin{equation*}
\begin{split}
E(F, v \finite)= -\limsup_k & \sum_v \log|N(v)| \Big(
\frac{\sum_{i,v} x_{v,i,k} [K_{v,i,k}:K_v]}{d^{nk}}
\\ & -\deg(F)
\frac{y_{v,i,k} [K_{v,i,k}:K_v]}{d^{nk}}- v(F)\Big).
\end{split}
\end{equation*}
\end{remark}
\begin{remark}
In case we are working in dimension one, that is with a map
$\vfi=(p_0:p_1): \Pl^1 \rightarrow \Pl^1 $, we can do a more precise
computation of the contributions of $y_{i,v,k}$ (See lemma 5.9 in
\cite{sz-t-p}). In this case, assuming that $p_0=A_dT^d_0+...$, the
integral at a finite place $v$ takes the particular form
\begin{equation*}
\begin{split}
\int_{\C_v^n} \log|F|_v d\mu_{v,\varphi}= -\lim_k & \log|N(v)| \Big(
\lim_{k} \frac{x_{v,i,k} [K_{v,i,k}:K_v]}{d^k}
\\ & -\deg(F) \frac{A_{d}[K_{v,i,k}:K_v]}{d-1}- v(F)\Big)
\end{split}
\end{equation*}
and again
$$E(F, v \finite)=\sum_v \int_{\Pl^n_{\bar{K}_v^n}} \log|F|
d\mu_{v,\varphi}.$$
\end{remark}
\begin{remark}
The geometry of $E(F, v \finite)$ allows to express in a better way
the negativity conditions for a model, indeed if our model
$(p_0:...:p_n)$ have negativity conditions and we are able to pick
the equation $F$ such that $v(F)=0$ for every $v$, then $E(F, v
\finite) \leq 0$.
\end{remark}
\section{Mahler formula}
In this section we present the main result of this paper.  Suppose
that $\vfi : \Pl_K^n \rightarrow \Pl_K^n$ is a map on the
n-dimensional projective space, given by a model $(p_0:...:p_n) :
\Pl_K^n \rightarrow \Pl_K^n$, with the property that $(p_0,...,p_n)$
is a regular sequence in $\cO_K$.  Recall the non-standard models
$\sigma_k : X_k \rightarrow \Pl^n$
  that we introduced before as the blows up of $\Pl^n$ at the subschemes
$Y_k$.
  Denote by $K_{v,i,k}$ the field of definition of $C_{v,i,k}$ and by
$K_v$ the local field of $K$ at $v$. Recall that $F_k=\sigma^*_k F$
for a polynomial $F=F(T_0/T_n,...,T_{n-1}/T_n)$ and $\L_k=\vfi_k^*
\L_0$, for $\L_0=\matO(1)$. The divisors $\Div(F)=D-\deg(F) \infty +
\sum_{v \finite} v(F) X_v$ and $\Div(T_n)=\infty$ represent divisors
in the arithmetic variety $\Pl_{\cO_K}^n$.
\begin{theorem} \label{Mahler formula} With the conditions and notations of
the above paragraph, we have the equality: $$h_{\vfi}(D)=\sum_{v / \infty}
\int_{\Pl^n_{\C}} \log |F|_v
d\mu_{\varphi,v} + E(F, v \finite) + \deg(F) h_{\vfi}(\infty).$$
\end{theorem}
\begin{proof}
We are going to make use of the arithmetic intersection theory on
$X_k$. Let's compute $\hat{\deg}(\hat{c}_1(L_k)^n |\Div(F_k))$. We
have the following:
\begin{enumerate}
\item{$\hat{\deg}(\hat{c}_1(L_k)^n |\Div(F_k))=d^{nk}
\sum_{v|\infty}\int_{\mathbb{P}^{1}(\bC_v)}
  \log|F|_v c_1(v,\|.\|_k)^n$ by proposition \ref{inters-at-infinity}.}
\item{$ \hat{\deg}
(\hat{c}_1(L_k))^n | C_k ) = [K_{v,i,k}:K_v] \log|N(v)|$  by
proposition \ref{general blow-up}.}
\item{$ \hat{\deg}
(\hat{c}_1(L_k))^n | X_{k,v} ) = \log|N(v)|$  also by proposition
\ref{general blow-up}.}
\end{enumerate}
Let's recall the formula \ref{Mahler-Step-k}:
\begin{equation*}
\begin{split}
\Div(F_k) =  D_k  - \deg(F) \infty_k & + \sum_{v,i} x_{v,i,k}
C_{v,i,k}  \\ & -\deg(F) \sum_{v,i} y_{v,i,k} C_{v,i,k} +
\sum_{\finite \,v} v(F) X_{v,k}.
\end{split}
\end{equation*}
Now we are going to let $\hat{\deg}(\hat{c}_1(L_k)^n |.)$ act on
each side,
\begin{equation*}
\begin{split}
\hat{\deg}(\hat{c}_1 (L_k)^n |\Div(F_k))  & = h_{\L^n_k}(D_k) -
(\deg(F)) h_{\L^n_k}(\infty_k) \\ &  + \sum_{i,v} x_{v,i,k} \log
|N(v)| [K_{v,i,k}:K_v] \\ & - \deg(F) \sum_{i,v} y_{v,i,k}  \log
|N(v)| [K_{v,i,k}:K_v] \\ & + d^{nk} \sum_{v} v(F) \log|N(v)|,
\end{split}
\end{equation*}
dividing by $d^{nk}$ and taking limits gives us that the limit $$\lim_k
\sum_{i,v} (x_{v,i,k} \log
|N(v)| [K_{v,i,k}:K_v] - \deg(F) y_{v,i,k} \log
|N(v)|[K_{v,i,k}:K_v])$$ exists and
\begin{equation*}
\sum_{v|\infty}\int_{\mathbb{P}^{n}(\bC)}
  \log|F|_v d\mu_{\vfi,v} =  - E(F, v \finite)  + h_{\vfi}(D_{\Q})-
  \deg(F)h_{\vfi}(\infty),
\end{equation*}
which was the result we wanted to prove.
\end{proof}
\begin{corollary}
Let $F$ be a rational function on $\Pl^n_K$ and $\Div(F)=D^+-D^-$
then
$$h_{\vfi}(D^+)-h_{\vfi}(D^-)=\sum_{v/\infty} \int_{\Pl^n_{\C}} \log
|F|_v d\mu_{\varphi,v} + E(F,v \finite) $$
\end{corollary}
\begin{proof}
The rational function $F$ can be written as the quotient $F^+/F^-$
of two homogeneous polynomial equations $F^+$ and $F^-$ of the same
degree.  Then we apply the previous result to $F^+$ and $F^-$.
\end{proof}
\begin{corollary}
If $\vfi$ has a model such that the divisor $\infty$ has a finite
forward orbit $\{\infty, \vfi(\infty),... \}$ (which forces
$h_{\vfi}(\infty)=0$), then
$$h_{\vfi}(D_{\Q})=\sum_{v / \infty} \int_{\Pl^n_{\C}} \log
|F|_v d\mu_{\varphi,v} + E(F, v \finite).$$
\end{corollary}
\begin{corollary}
Suppose that $\vfi=(p_0:p_1) : \Pl^1 \rightarrow \Pl^1$ and we have
chosen coordinates such that $T_n/p_1$. The integral of the log the
minimal equation $F$ of a point $P=(\lambda:1) \in \Pl^1$ is related
to the height of $P$ by the formula
$$h_{\vfi}(P)=\frac{1}{\deg(F)}\sum_v \int_{\Pl^1_{\bar{K_v}}} \log
|F|_v d\mu_{\varphi,v}.$$
\end{corollary}
\begin{proof}
This is a combination of the previous corollary, proposition
\ref{convergence of finite place int} and definition \ref{definition
of int at finite place}.
\end{proof}
\begin{corollary}
Assume $h_{\vfi}(\infty)=0$.  Let $(p_0:...:p_n)$ be a model for
$\vfi : \Pl^n \rightarrow \Pl^n$, satisfying negativity conditions
and with $v(F)=0$ for every finite place $v \in \cO_k$.  Then
$$h_{\vfi}(D_{\Q}) \leq \sum_{v/\infty} \int_{\Pl^n_{\C}} \log
|F|_v d\mu_{\varphi,v}.$$
\end{corollary}
\begin{example}
The rational map $\vfi : \Pl_{\Q}^2 \rightarrow \Pl_{\Q}^2$ given by
the model $\vfi(x:y:z)=(x^2+yx:y^2+zx+zy:z^2)$, has good reduction
everywhere, that is, the expression $E(F,v \finite)=0$.
\end{example}
\begin{example}
The rational map $\vfi : \Pl_{\Q}^2 \rightarrow \Pl_{\Q}^2$ given by
the model $\vfi(x:y:z)=(y^2-3z^2:x^2-3y^2:zy)$ has bad reduction at
$3$.  The reduced map is not defined at the point $(0:0:1) \in
\Pl^2_{\bar{\Q}_3}$. The model satisfy negativity conditions and
$h_{\vfi}(\infty)=0$. A polynomial equation $F=c_m z^m + ...$ has
$c_m \equiv 0 (mod$ $3)$ if and only if $E(F,v \finite) \neq 0$.
\end{example}
\begin{example}
The rational map $\vfi : \Pl_{\Q}^2 \rightarrow \Pl_{\Q}^2$ given by
the model $\vfi(x:y:z)=(3y^2-5z^2:3x^2-5y^2:zy)$ has bad reduction
at the places $3$ and $5$.  The reduced maps is not well defined at
the points $(1:0:0) \in \Pl^2_{\bar{\Q}_3}$ and $(0:0:1) \in
\Pl^2_{\bar{\Q}_5}$.
\end{example}
\begin{example}
Let $p$ be a prime number.  The map $\vfi : \Pl^1_{\Q} \rightarrow
\Pl^1_{\Q}$ given by the model $\vfi(x:y)=(px^2+y^2:py^2)$ has bad
reduction at $p$.  Following lemma 5.9 in \cite{sz-t-p}, if $F$ is a
polynomial equation, we have $E(F,v \finite)= \log(p)$.

\end{example}

\end{document}